\documentclass[a4paper,12pt]{article}
    \usepackage[top=2.5cm,bottom=2.5cm,left=2.5cm,right=2.5cm]{geometry}
    \usepackage{cite, amsmath, amssymb}
\begin{document}
\begin{center}
{\LARGE\bf $\mathbf{L}$-Borderenergetic Graphs}
\end{center}
$\vspace{0,5cm}$
\begin{center}
{\large \bf Fernando Tura}
\end{center}
\begin{center}
\it Departamento de Matem\'atica, UFSM, Santa Maria, RS, 97105-900, Brazil
\end{center}
\begin{center}
\tt ftura@smail.ufsm.br
\end{center}

\newtheorem{Thr}{Theorem}
\newtheorem{Pro}{Proposition}
\newtheorem{Que}{Question}
\newtheorem{Con}{Conjecture}
\newtheorem{Cor}{Corollary}
\newtheorem{Lem}{Lemma}
\newtheorem{Fac}{Fact}
\newtheorem{Ex}{Example}
\newtheorem{Def}{Definition}
\newtheorem{Prop}{Proposition}
\def\floor#1{\left\lfloor{#1}\right\rfloor}

\newenvironment{my_enumerate}{
\begin{enumerate}
  \setlength{\baselineskip}{14pt}
  \setlength{\parskip}{0pt}
  \setlength{\parsep}{0pt}}{\end{enumerate}
}

\newenvironment{my_description}{
\begin{description}
  \setlength{\baselineskip}{14pt}
  \setlength{\parskip}{0pt}
  \setlength{\parsep}{0pt}}{\end{description}
}


\begin{abstract}
The  energy of a graph is defined as the sum the absolute values of the eigenvalues  of its adjacency matrix.
A graph $G$ on $n$ vertices is said to be borderenergetic  if its energy equals the energy of the complete graph $K_n.$
In this paper, we promote this concept for the Laplacian matrix.
The Laplacian energy of $G,$ introduced by Gutman and Zhou \cite{Gutman}, is given by  $LE(G)= \sum_{i=1}^n | \mu_i - \overline{d} |,$ where $ \mu_i$ are the Laplacian eigenvalues of $G$ and $\overline{d}$ is the average degree of $G.$ In this way, we say $G$ to be  $L$-borderenergetic if $LE(G)= LE(K_n).$
Several  classes of $L$-borderenergetic  graphs are obtained including  result that for each  integer $r\geq 1,$ there are $2r+1$ graphs, of order $n=4r+4,$ pairwise  $L$-noncospectral and  $L$-bordernergetic graphs.
\end{abstract}

\baselineskip=0.30in

\begin{center}
     (Received April 2, 2016)
\end{center}

\section{Introduction}
\label{intro}
Throughout this paper, all graphs are assumed to be finite, undirected and without loops or multiple edges.
If $G$ is a graph of order $n$ and $M$ is a real symmetric matrix associated with $G,$ then  the $M$- {\em energy} of $G$ is
\begin{equation}
\label{energy}
 E_{M}(G)= \sum_{i=1}^n |  \lambda_i (M) - \frac{ tr(M)}{n}|.
 \end{equation}
The  energy of a graph  simply refers to using  the adjacency matrix in (\ref{energy}).
There are many results on energy \cite{Li, Li2, Li3, Li4, Li5, Li6} and its applications in several areas, including in chemistral
 see  \cite{Gutman2012} for more details and the references therein.

It is  well known that the complete graph $K_n$ has  $E(K_n) = 2n-2.$ In this context, several authors have been presented families of graphs
with same energy of the complete graph $K_n.$ Recently, Gong, Li,  Xu, Gutman and Furtula   \cite{Gutman2015} introduced the concept of {\em bordernergetic}. A graph $G$ on $n$ vertices is said to be borderenergetic  if its energy equals the energy of the complete graph $K_n.$

In \cite{Gutman2015}, it was shown that there exits borderenergetic graphs  on order $n$ for each integer $n\geq 7,$ and all borderenergetic graphs with $7,8,$ and $9$ vertices were determined.

In \cite{JTT2015} considered the eigenvalues and energies of threshold graphs. For each $n\geq 3,$ they determined $n-1$ threshold graphs on $n^2$ vertices, pairwise non-cospectral and equienergetic to the complete graph
$K_{n^2}.$
Recently,  Hou and Tao  \cite{Hou}, showed that for each $n\geq 2$ and $p\geq 1$  $(p \geq 2$ if $n=2),$ there are  $n-1$ threshold graphs on $pn^2$ vertices, pairwise non-cospectral and equienergetic  with the complete graph
$K_{pn^2},$
generalizing the results in \cite{JTT2015}.

The Laplacian energy of $G,$ introduced by Gutman and Zhou \cite{Gutman}, is given by
\begin{equation}
LE(G)= \sum_{i=1}^n | \mu_i - \overline{d} |
\end{equation}
where $ \mu_i$ are the Laplacian eigenvalues of $G$ and $\overline{d}$ is the average degree of $G.$ Similarly for the laplacian energy, we have that $LE(K_n)= 2n-2.$

The first purpose of this paper is to promote the concept of borderenergetic to the laplacian matrix.
In this way, we say $G$ to be  $L$-borderenergetic if $LE(G)= LE(K_n).$
The second is to present  several  classes of  $L$-borderenergetic graphs.

The paper is organized as follows.
In Section 2 we describe some  known results about the Laplacian spectrum of graphs.
In Section 3 we present four classes of  $L$-borderenergetic. We finalize this paper, showing that for each  integer $r\geq 1,$ there are $2r+1$ graphs, of order $n=4r+4,$ pairwise  $L$-noncospectral and  $L$-bordernergetic graphs.

\section{Premilinares}
\label{Diag}

Let $G_1= (V_1, E_1)$ and $G_2= (V_2, E_2)$  be undirected graphs without  loops or multiple edges.
The {\em union} $G_1 \cup G_2$ of graphs $G_1$ and $G_2$  is the graph $G=(V,E)$  for which $V= V_1 \cup V_2$ and $E = E_1 \cup E_2.$
We denote the graph  $ \underbrace{G  \cup G \cup \ldots \cup G}_{m} $ by  $m G.$
The {\em join} $G_1 \nabla G_2$  of graphs $G_1$ and $G_2$  is the graph obtained from $G_1 \cup G_2$ by joining every vertex of $G_1$
with every vertex of $G_2.$

The Laplacian spectrum of $G_1  \cup \ldots \cup G_k$ is the union of Laplacian spectra of $G_1, \ldots, G_k,$ while
the Laplacian spectra of the complement  of $n$- vertex graph $G$ consists  of values $n - \mu_i,$ for each Laplacian eigenvalue
$\mu_i$ of $G,$ except for a single instance of eigenvalue $0$ of $G.$


\begin{Lem}
\label{theorem1}
Let $G$ be a graph  on  $n$  vertices with  Laplacian matrix  $L.$  Let  $0=\mu_1 \leq \mu_2 \leq \ldots \leq \mu_{n} $ be
 the eigenvalues of $L.$ Then the eigenvalues of $ \overline{G}$ are
$$ 0 \leq n-\mu_n  \leq n - \mu_{n-1} \leq  n - \mu_{n-2} \leq \ldots \leq  n - \mu_{2}$$
with the same corresponding  eigenvectors.
\end{Lem}
{\bf Proof:}
Note that the Laplacian matrix of $\overline{G}$ satisfies $L(\overline{G})=  nI + J - L,$  where $I$ is the identity matrix and $J$ is the matrix each of whose
entries is equal 1. Therefore, for $i=2, \ldots, n,$ if $x$ is an eigenvector of $L$ corresponding to $\mu_i,$ then $Jx=0.$ Therefore
$$ L(\overline{G})x =  (nI + J - L)x =  nIx + Jx - Lx = (n - \mu_i) x.$$
Thus $ n - \mu_i$ is an eigenvalue with $x_i$ as a corresponding eigenvector. Finally, $e=(1,\ldots,1)$ is an eigenvector of $L(\overline{G})$
corresponding to $0.$ $\hspace{7,5cm} \square$

Recall that $G$ is laplacian integral if its spectrum consists entirely of integers \cite{Kirkland, Merris2}. Follows  from Lemma \ref{theorem1} that $G$ is laplacian integral  if and only if  $\overline{G}$ is  laplacian integral.

\begin{Thr}
\label{theorem2}
Let $G_1$ and $G_2$ be graphs on  $n_1$ and $n_2$ vertices, respectively. Let $L_1$ and $L_2$ be the Laplacian matrices for $G_1$ and $G_2,$
respectively, and let $L$ be the Laplacian matrix for $G_1 \nabla G_2.$  If $0=\alpha_1 \leq \alpha_2 \leq \ldots \leq \alpha_{n_1} $ and $0=\beta_1 \leq \beta_2 \leq \ldots \leq \beta_{n_2}$
are the eigenvalues of $L_1$ and $L_2,$ respectively. Then the eigenvalues of $L$ are
$$ 0, \hspace{0,2cm} n_2 + \alpha_2, \hspace{0,2cm} n_2 + \alpha_3 , \ldots, \hspace{0,1cm} n_2 + \alpha_{n_1}$$
$$ n_1 + \beta_2, \hspace{0,2cm} n_1 + \beta_3, \ldots, \hspace{0,1cm} n_1 + \beta_{n_2}, \hspace{0,2cm} n_1 + n_2.$$
\end{Thr}
{\bf Proof:}
Since that the join of graphs $G_1$ and $G_2$ is given by $G_1\nabla G_2=   \overline{  \overline{G_1} \cup \overline{G_2}}$ (see \cite{Gutman}), the proof follows immediately from the Lemma \ref{theorem1}. $\hspace{7,5cm} \square$

\section{$\mathbf{L}$-Borderenergetic graphs}
\label{SecInterval}
Recall that the $ L$-energy   of a graph  $G$ is  obtained by
$LE(G)= \sum_{i=1}^n | \mu_i - \overline{d} |,$
 where $  \mu_i$  are the laplacian eigenvalues  of $G$ and $\overline{d}$ is the average degree of $G.$
 It is known that the complete graph $K_n$  has Laplacian energy  $2(n-1).$
We exhibit four infinite classes  $\Omega_i =\{G_1, G_2, \ldots, G_r, \ldots \}$ for $i=1,\ldots,4$ such that each $G_r,$ of order $n=4r+4,$ satisfies $LE(G_r)= LE(K_{4r+4}).$

\subsection{ The class $\Omega_1$ }
For  each integer $r\geq 1,$  we define the graph $G_r \in \Omega_1$ to be the following join
$$ G_r = (rK_1 \cup ( K_1 \nabla (r+1) K_1)) \nabla (rK_1 \cup ( K_1 \nabla (r+1) K_1)).$$
$G_r$ has order $n=4r+4.$ We let $\mu^m$ denote  the laplacian eigenvalue $\mu $ with multiplicity equals to $m.$

\begin{Lem}
\label{lema1}
Let $G_r \in \Omega_1$ be a graph of order $n =4r+4.$
Then  the Laplacian spectrum of $ G_r$ is given by
$$ 0 ;  \hspace{0,2cm} (2r +2)^{2r} ;\hspace{0,2cm} (2r+3)^{2r}; \hspace{0,2cm} (3r+4)^2; \hspace{0,2cm}  4r +4. $$

\end{Lem}
{\bf Proof:}
Let $G_r \in \Omega_1.$ Let's denote $ H = rK_1 \cup (K_1 \nabla (r+1)K_1).$ By definition we have that
$ G_r =  H \nabla  H.$
According by Theorem \ref{theorem2}, we just need to determine the Laplacian spectrum of the $H$ and add its order.
By direct calculus follows that the Laplacian spectrum of $H$ is equal to
 $$ 0^r ;  \hspace{0,2cm} 1^{r} ;\hspace{0,2cm} r+2.$$
Since $H$ has order $2r+2,$ by Theorem \ref{theorem2} the result follows. $\hspace{5,5cm}\square$

\begin{Thr}
\label{}
For each $r\geq 1,$  $G_r$ is  $L$-borderenergetic and  L-noncospectral graph with $K_{4r+4}.$
\end{Thr}
{\bf Proof:}
Clearly $G_r$ and $K_{4r+4}$ are  L-noncospectral.
Let $\overline{d}$ be the average degree of $G_r.$ Since that  $\overline{d}$ is equal to average of Laplacian eigenvalues of $G_r$ then
 $ \overline{d} = \frac{ 2r(4r+5) + 2(3r+4) +4r+4}{4r+4}\\ = 2r+3.$ Using Lemma \ref{lema1},
 $ LE(G_r) = 4r+4 -(2r+3) +2(3r+4 -2r-3) +2r(2r+3 -2r-3)+2r(2r+3- 2r-2) + 2r+3 = 8r+6 = LE(K_{4r+4}). $ $\hspace{7,25cm} \square$

\subsection{ The class $\Omega_2$ }
For  each integer $r\geq 1,$  we define the graph $G_r \in \Omega_2$ to be the following join
$$G_r =  (r+1) K_2 \nabla (r+1) K_2.$$
$G_r$ has order $n=4r+4.$ We let $\mu^m$ denote  the laplacian eigenvalue $\mu $ with multiplicity equals to $m.$

\begin{Lem}
\label{lema2}
Let $G_r \in \Omega_2$ be a graph of order $n =4r+4.$
Then  the Laplacian spectrum of $ G_r$ is given by
$$ 0 ;  \hspace{0,2cm} (2r +2)^{2r} ;\hspace{0,2cm} (2r+4)^{2r +2}; \hspace{0,2cm} 4r +4.$$

\end{Lem}
{\bf Proof:}
Let $G_r \in \Omega_2.$ Let's denote $ H = ( r+1)K_2 .$ By definition we have that
$ G_r =  H \nabla  H.$
According by Theorem \ref{theorem2}, we just need to determine the Laplacian spectrum of the $H$ and add its order.
By direct calculus follows that the Laplacian spectrum of $H$ is equal to
 $$ 0^{r+1} ;  \hspace{0,2cm}  2^{r+1} .$$
Since $H$ has order $2r+2,$ by Theorem \ref{theorem2} the result follows. $\hspace{5,5cm} \square$

\begin{Thr}
\label{}
For each $r\geq 1,$  $G_r$  is  $L$-borderenergetic and  L-noncospectral graph with $K_{4r+4}.$
\end{Thr}
{\bf Proof:}
Clearly $G_r$ and $K_{4r+4}$ are  L-noncospectral.
Let $\overline{d}$ be the average degree of $G_r.$ Since that  $\overline{d}$ is equal to average of Laplacian eigenvalues of $G_r$ then
 $ \overline{d} = \frac{ (2r+2)(4r +4)+4r+4}{4r+4} = 2r+3.$ Using Lemma \ref{lema1},
 $ LE(G_r) = 4r+4 -(2r+3) +(2r +2)(2r+4 -2r-3) +2r(2r+3 -2r-2)+ 2r+3 = 8r+6 = LE(K_{4r+4}). $
$\hspace{7,5cm} \square$

\subsection{ The classes $\Omega_3$ and $\Omega_4$ }
For  each integer $r\geq 1,$  we define the following two graphs $G_r \in \Omega_3$ and $G'_r \in \Omega_4:$
$$G_r = (K_2 \cup (2r+1)K_1) \nabla (2r+1)K_1,$$
$$G'_r = ((2r+1)K_1) \nabla (2r+2)K_1  \nabla K_1,$$
where $G_r$ and $G'_r$ have order $n=4r+4.$

The proof of following results are similar to others above, then we will omite them.

\begin{Lem}
\label{lema3}
Let $G_r \in \Omega_3$ and $G'_r \in \Omega_4$ be graphs of order $n =4r+4.$
Then  the Laplacian spectrum of $ G_r$ and $G'_r$ are given by
$$ 0 ;  \hspace{0,2cm} (2r +1)^{2r+1} ;\hspace{0,2cm} (2r+3)^{2r +1}; \hspace{0,2cm} 4r +4,$$
$$ 0 ;  \hspace{0,2cm} (2r +2)^{2r+1} ;\hspace{0,2cm} (2r+3)^{2r}; \hspace{0,2cm} (4r +4)^2,$$
respectively.
\end{Lem}

\begin{Thr}
\label{}
For each $r\geq 1,$  $G_r$  and $G'_r$   are  $L$-borderenergetic and  L-noncospectral graphs with $K_{4r+4}.$
\end{Thr}

\section{More $\mathbf{L}$-Borderenergetic graphs}
\label{secEquienergetic}

In this Section we obtain more $L$-borderenergetic graphs including  result that for each integer $r\geq 1,$ there are $2r+1$ graphs, of order $n=4r+4,$ pairwise  $L$-noncospectral and  $L$-bordernergetic graphs.
Consider the following graphs:
$$ H_1 = rK_1 \cup ( K_1 \nabla (r+1) K_1)$$
$$ H_2 = (r+1) K_2 $$
$$ H_3 =  r K_2 \cup 2K_1 $$
$$ H_4 =  ((2r+1)K_1)  \nabla K_1.$$
The proof of following results are similar to others above, then we will omite them.

\begin{Lem}
\label{}
Let $G_{1,2}$ be a graph of order $n =4r+4$ obtained by the following join $G_{1,2} = H_1 \nabla H_2.$
Then  the Laplacian spectrum of $ G_{1,2}$ is given by
$$ 0 ;  \hspace{0,2cm} (2r +2)^{2r} ;\hspace{0,2cm} (2r+3)^{r}; \hspace{0,2cm} (2r+4)^{r+1};\hspace{0,2cm} 3r +4; \hspace{0,2cm} 4r +4.$$
\end{Lem}

\begin{Lem}
\label{}
Let $G_{1,3}$ be a graph of order $n =4r+4$ obtained by the following join $G_{1,3} = H_1 \nabla H_3.$
Then  the Laplacian spectrum of $ G_{1,3}$ is given by
$$ 0 ;  \hspace{0,2cm} (2r +2)^{2r +1} ;\hspace{0,2cm} (2r+3)^{r}; \hspace{0,2cm} (2r+4)^{r};\hspace{0,2cm} 3r +4; \hspace{0,2cm} 4r +4.$$
\end{Lem}

\begin{Lem}
\label{}
Let $G_{2,3}$ be a graph of order $n =4r+4$ obtained by the following join $G_{2,3} = H_2 \nabla H_3.$
Then  the Laplacian spectrum of $ G_{2,3}$ is given by
$$ 0 ;  \hspace{0,2cm} (2r +2)^{2r +1} ; \hspace{0,2cm} (2r+4)^{2r +1}; \hspace{0,2cm} 4r +4.$$
\end{Lem}

\begin{Lem}
\label{}
Let $G_{2,4}$ be a graph of order $n =4r+4$ obtained by the following join $G_{2,4} = H_2 \nabla H_4.$
Then  the Laplacian spectrum of $ G_{2,4}$ is given by
$$ 0 ;  \hspace{0,2cm} (2r +2)^{r} ;\hspace{0,2cm} (2r+3)^{2r}; \hspace{0,2cm} (2r+4)^{r +1}; \hspace{0,2cm} (4r +4)^2.$$
\end{Lem}

\begin{Lem}
\label{}
Let $G_{3,4}$ be a graph of order $n =4r+4$ obtained by the following join $G_{3,4} = H_3 \nabla H_4.$
Then  the Laplacian spectrum of $ G_{3,4}$ is given by
$$ 0 ;  \hspace{0,2cm} (2r +2)^{r+1} ;\hspace{0,2cm} (2r+3)^{2r}; \hspace{0,2cm} (2r+4)^{r}; \hspace{0,2cm} (4r +4)^2.$$
\end{Lem}

\begin{Thr}
\label{}
For each integer $r\geq 1,$  $G_{1,2}, G_{1,3}, G_{2,3}, G_{2,4}$ and $G_{3,4}$  are  $L$-borderenergetic and  L-noncospectral  graphs.
\end{Thr}

For integers $r\geq 1$ and $i=0,1,\ldots, 2r,$  consider the following $2r+1$ graphs:
$$ G_{i,r} =  ((2r+1)K_1) \nabla ( (2r+1-i)K_1 )\cup  (K_1 \nabla (i+1) K_1), $$
of order $n=4r+4.$

\begin{Lem}
\label{}For integers $r\geq 1$ and $i=0,1,\ldots, 2r,$  let $G_{i,r}$ be a graph of order $n =4r+4.$
Then  the Laplacian spectrum of $ G_{i,r}$ is given by
$$ 0 ;  \hspace{0,2cm} (2r +1)^{2r +1-i} ;\hspace{0,2cm} (2r+2)^i; \hspace{0,2cm} (2r+3)^{2r}; \hspace{0,2cm} (2r+3+i); \hspace{0,2cm} (4r +4).$$
\end{Lem}

\begin{Thr}
\label{}For integers $r\geq 1$ and $i=0,1,\ldots,2r,$  $G_{i,r}$  are  $L$-borderenergetic and  $L$-noncospectral  graphs.
\end{Thr}

\end{document}